\documentstyle[12pt]{article}
\def\ba{\begin{array}}
\def\ea{\end{array}}	
\def\be{\begin{equation}}
\def\ee{\end{equation}}
\def\ot{\otimes}

\def\d{\delta}

\def\ep{\epsilon}

\def\si{\sigma}

\def\th{\Theta}

\def\Om{\Omega}
\def\ra{\longrightarrow}

\def\ca{{\cal A}}

\def\ca{{\cal A}}

\def\cm{{\cal M}}

\def\cw{{\cal W}}

\def\dn#1,{_{({#1})}}
\newfont{\numb}{msbm10}
\def\com{\mbox{\numb C}}

\def\pg#1,#2,#3,{\langle #1 | #2 \rangle^{#3}}
\def\otn#1,#2,{\left( #1 \otimes #2 \right) }

\def\ps#1,#2,{\Psi_{#1,#2}}
\def\id#1,{id_{#1}}
\def\ev#1,{ev_{#1}}
\def\tl#1,{\tilde{#1}}
\def\1n{^{(1)}}
\def\2n{^{(2)}}
\def\3{^{(3)}}
\def\ps#1,#2,{\Psi_{{#1}{,}{#2}}}
\def\et#1,{U^{\ot #1}}
\def\get#1,{{\cal U}^{\ast \ot #1}}
\def\pr{(.|.)}
\begin{document} \title{Remarks on Quantum Statistics} 
\author{Wladyslaw Marcinek
\thanks{The work is partially sponsored 
by Polish Committee for Scientific Research (KBN) under 
Grant 2P03B130.12}
\\ 
Institute of Theoretical Physics, University
of Wroc{\l}aw, Poland} 
\maketitle 
\begin{abstract} 
Some problems related to an algebraic approach to quantum 
statistics are discussed. Quantum statistics is described
as a result of interactions. The Fock space representation
is discussed. The problem of existence of well--defined 
scalar product is considered. An example of physical 
applications is also given.
\end{abstract} 
\section*{Introduction} 
In the last years a few different approaches to quantum 
statistics which generalize the usual boson or fermion 
statistics has been intensively developed by several
authors. The so--called $q$--statistics and corresponding 
$q$--relations have been studied by Greenberg \cite{owg,gre}, 
Mohapatra \cite{moh}, Fivel \cite{fi} and many others, see
\cite{zag,mepe,bs2} for example. The deformation of 
commutation relations for bosons and fermions corresponding 
to quantum groups $SU_q(2)$ has been given by Pusz and 
Woronowicz \cite{P,PW}. The $q$--relations corresponding 
to superparticles has been considered by Chaichian, 
Kulisch and Lukierski \cite{ckl}. Quantum deformations 
have been also studied by Vokos \cite{V}, 
Fairle and Zachos \cite{FZ} and many others.

Note that there is also an approach to particle systems 
with some nonstandard statistics in low dimensional spaces 
based on the notion of the braid group $B_n$ \cite{Wu,I}.
In this approach the configuration space for the system of 
$n$--identical particles moving on a manifold $\cm$ is given 
by the formula
$$
Q_{n}(\cm) = {\left(\cm^{\times n} - D\right) }/{S_{n}},
$$ 
where $D$ is the subcomplex of the Cartesian product 
$\cm^{\times n}$ on which two or more particles occupy 
the same position and $S_{n}$ is the symmetric group. 
The group $\pi_{1}\left( Q_{n}(M)\right)\equiv B_{n} (M)$
is known as the $n$--string braid group on $\cm$. Note that 
the group $\Sigma_n (\cm)$ is a subgroup of $B_n(\cm)$ and 
is an extension of the symmetric group $S_n$ describing the 
interchange process of two arbitrary indistinguishable 
particles. It is obvious that the statistics of the given 
system of particles is determined by the group $\Sigma_n$ 
\cite{Wu,I}. The mathematical formalism related to the braid 
group statistics has been developed intensively by Majid, see
\cite{SM,Maj,SMa,bm,sma,qm,ssm} for example. 
It is interesting that all commutation relations for particles 
equipped with arbitrary statistics can be described as 
representations of the so-called quantum Weyl or Wick algebra 
$\cw$. Such algebraic formalism has been considered by J\"orgensen,
Schmith and Werner \cite{jswe} and further developed by
the author in the series of papers 
\cite{WM3,WM4,WM6,WM8,wm7,mad,mco,qweyl,wmq}. 
and also by Ralowski \cite{jswe,RM,m10,ral}.
Similar approach has been also considered 
by others authors, see \cite{sci,twy,mep} and \cite{mphi,mira}.
An interesting approach to quantum statistics has been also given 
in \cite{fios,melme}. A proposal for the general algebraic formalism 
for description of particle systems equipped with an arbitrary 
generalized statistics based on the concept of monoidal categories 
with duality has been given by the author in  \cite{castat}.
The physical interpretation for this formalism was shortly 
indicated. A few examples of applications for this formalism 
are considered in \cite{top,sin}. 

The generalization of the concept of quantum statistics
is motivated by many different applications in quantum field 
theory and statistical physics. Some problems in condensed 
matter physics, magnetism or quantum optics lead to the
study of particle systems obeying nonstandard statistics.
It is interesting that in the last years new and highly 
organized structures of matter has been discovered. 
For example in fractional quantum Hall effect a system 
with well defined internal order has appears \cite{zee}. 
Another interesting structures appear in the so called 
$\frac{1}{2}$ electronic magnetotransport
anomaly \cite{jai,dst}, high temperature superconductors 
or laser excitations of electrons. In these cases
certain anomalous behaviour of electron have appear.
An another example is given by the so--called Lutinger 
liquid \cite{hal1}. The concept of statistical--spin 
liquids has been studied by Byczuk and Spalek \cite{bys}. 
It is interesting that half of the available 
single--particle states are removed by the statistical 
interaction between the particles with opposite spins. 
The study of highly organized structures leads to 
the investigation of correlated systems
of interacting particles. The essential problem in 
such study is to transform the system of interacting
particles into an effective model convenient for
the description of ordered structures. 
A system with generalized statistics seems to be
one of the best candidate for such model.
Hence there is an interest in the development of 
formalism related to the particle systems
with unusual statistics and the possible
physical applications.
In this paper we would like to discuss some problems
of possibility of application of systems with generalized 
statistics for the further description of ordered structures.
All our considerations are based on the assumption that 
the quantum statistics of charged particles is determined 
by some specific interactions.
\section*{Fundamental Assumptions} 
The starting point for our discussion 
is a system of charged particles interacting with certain 
quantum field. The proper physical nature of the system is 
not essential for our considerations. 
The fundamental assumption is that the problem of 
interacting particles can be reduced to the study 
of a system consisting $n$ charged particles and 
$N$--species of quanta of the field. In this way we 
can restrict our attention to study of such system. 
It is natural to expect that some new excited states 
of the system have appear as a result of certain specific
interaction. The existence of new ordered structures
depends on the existence of such additional excitations.
Hence we can restrict our attention to the study
of possibility of appearance for these excitations.
For the description of such possible excited states 
we use the concept of dressed particles.
We assume that every charged particle is equipped with 
ability to absorb quanta of the external field. A system
which contains a particle and certain number of quanta as
a result of interaction with the external field is said to 
be dressed particle. A particle without quantum is called 
undressed or a quasihole. The particle dressed with two quanta 
of certain species is understand as a system of two new objects
called quasiparticles. A quasiparticle is in fact the charged 
particle dressed with a single quantum. Two quasiparticles
are said to be identical if they are dressed with quanta
of the same species. In the opposite case when the particle 
is equipped with two different species of quanta then we have 
different quasiparticles. We describe excited states as 
composition of quasiparticles and quasiholes. 
It is interesting that quasiparticles and quasiholes have 
also their own statistics. 
We give the following assumption for the algebraic
description of excitation spectrum of single dressed
particle.

\paragraph*{Assumption 0: The ground state.}
There is a state $|0> = {\bf 1}$ called the ground one.
There is also the conjugate ground state $<0| \equiv 1^{\ast}$.
This is the state of the system before intersection.

\paragraph*{Assumption 1: Elementary states.}
There is an ordered (finite) set of single quasiparticle states
\be 
S := \{x^i : i = 1,\ldots, N<\infty\}.
\ee
These states are said to be elementary (simple). They
represent elementary excitations of the system.
We assume that the set $S$ of elementary states forms 
a basis for a finite linear space $E$ over a field
of complex numbers $\com$. 

\paragraph*{Assumption 2: Elementary conjugate states.}
There is also a corresponding set of single quasihole states 
\be
S^{\ast} := \{x^{\ast i} : i = N, N-1,\ldots , 1\}.
\ee
These states are said to be conjugated. The set $S^{\ast}$ of 
conjugate states forms a basis for the complex conjugate space 
$E^{\ast}$. The pairing $\pr : E^{\ast} \ot E \ra \mbox{\numb C}$
is given by
\be
(x^{\ast i}|x^j) := \d^{ij}.
\ee

\paragraph*{Assumption 3: Composite states.}
There is a set of projectors
\be
\Pi_n : E^{\ot n} \ra E^{\ot n}
\ee
such that we have a $n$--multilinear mapping
\be
\odot_n : E^{\times n} \ra E^{\ot n}.
\ee
defined by the following formula
\be
x^{i_1} \odot \cdots \odot x^{i_n} := 
\Pi_{n}(x^{i_1} \ot \cdots \ot x^{i_n}).
\ee
The set of $n$--multiquasiparticle states is denoted by
$P^n (S)$. All such states are result of composition
(or clustering) of elementary ones. These states are 
also called composite states of order $n$. They represent
additional excitations charged particle under interaction.
In this way for multiquasiparticle states we have the 
following set of states
\be
P^n (S) := \{x^{\si} \equiv x^{i_1} \odot \cdots \odot x^{i_n} 
: \si = (i_1 ,\ldots,i_n) \in I\},
\ee
Here $I$ is a set of sequences of indices such that the
above set of states forms a basis for a linear space
$\ca^n$. We have
\be
\ca^n = Im(\Pi_n).
\ee
Obviously we have where $\ca^0 \equiv {\bf 1}\com$, 
$\ca^1 \equiv E$ and $\ca^n \subset E^{\ot n}$.

\paragraph*{Assumption 4: Composite conjugated states.}
We also have a set of projectors
\be
\Pi^{\ast}_n : E^{\ast\ot n} \ra E^{\ast\ot n}
\ee
and the corresponding set of composite 
conjugated states of length $n$
\be
P^{n}(S^{\ast}) := \{x^{\ast \si} \equiv x^{\ast i_n} 
\odot \cdots \odot x^{\ast i_1} : 
\si = (i_1 ,\ldots,i_n) \in I \}.
\ee
The set $P^n (S^{\ast})$ of composite conjugated states of 
length $n$ forms a basis for a linear space $\ca^{\ast n}$.

\paragraph*{Assumption 5: Algebra of states.}
The set of all composite states of arbitrary length is denoted 
by $P(S)$. For this set of states we have the following
linear space
\be
\ca := \bigoplus\limits_{n} \ \ca^n .
\ee
If the formula
\be
m(s \ot t)s \equiv s \odot t := \Pi_{m+n}(\tl s, \ot \tl t,)
\label{mul}
\ee
for $s=\Pi_m (\tl s,), t=\Pi_n (\tl t,)$, 
$\tl s, \in E^{\ot n}, \tl t, \in E^{\ot m}$,
defines an associative multiplication in $\ca$, then we say 
that we have an algebra of states. This algebra represents
excitation spectrum for single dressed particle.

\paragraph*{Assumption 6: Algebra of conjugated states.}
The set of composite conjugated states of arbitrary length 
is denoted $by P(S^{\ast})$. We have here a linear space
\be
\ca^{\ast} := \bigoplus\limits_{n} \ \ca^{\ast n} ,
\ee
If $m$ is the multiplication in $\ca$, then the
multiplication in $\ca^{\ast}$ corresponds to 
the opposite multiplication in $\ca$
\be
m^{op}(t^{\ast} \ot s^{\ast}) = (m(s \ot t))^{\ast}.
\ee
\section*{Creation and Annihilation Operators} 
We define creation operators for our model as multiplication
in the algebra $\ca$
\be
a^+_{s} t := s \odot t, \quad \mbox{for} \quad s, t \in \ca,
\ee
where the multiplication is given by $\ref{mul}$.
For the ground state and annihilation operators we assume that
\be
\langle 0|0 \rangle = 0, \quad a_{s^{\ast}} |0\rangle = 0 
\quad \mbox{for} \quad s^{\ast} \in \ca^{\ast}.
\ee
The proper definition of action of annihilation 
operators on the whole algebra $\ca$ is a problem.
For the pairing 
$\pg -, -, n, : \ca^{\ast n} \ot \ca^n \ra \com$
we assume in addition that we have the following
formulae
\be
\pg 0, 0, 0, := 0,\quad \pg i, j, 1, := 
(x^{\ast i}|x^j) = \d^{ij},\\
\pg s, t, n, := \pg \tilde{s}, P_n \tilde{t}, n, _0
\quad\mbox{for}\quad n\geq2
\ee
where $\tl s,  , \tl t, \in E^{\ot n}$, 
$P_n : E^{\ot n}\ra E^{\ot n}$ is an additional linear
operator and
\be
\langle i_1 \cdots i_n |j_1\cdots j_n \rangle^n_0 :=
\pg i_1, j_1, 1, \cdots \pg j_n, j_n, 1, .
\ee
Observe that we need two sets $\Pi := \{\Pi_n\}$ and 
$P := \{P_n\}$ of operators and the action 
\be
a: s^{\ast} \ot t \in \ca^{\ast k} \ot \ca^n \ra 
a_{s^{\ast}} t \in \ca^{n-k}.
\label{act}
\ee
of annihilation operators for the algebraic description of 
our system. In this way the triple $\{\Pi , P, a\}$, where 
$\Pi$ and $P$ are set of linear operators and $a$ is the action 
of annihilation operators, is the initial data for our model. 
The problem is to find and classify all triples of initial
data which lead to the well--defined models. The general 
solution for this problem is not known for us. Hence we 
must restrict our attention for some examples.

\paragraph*{Definition:}
If operators $P$ and $\Pi$ and the action $a$ of
annihilation operators are given in such a way that 
there is unique, nondegenerate, positive definite 
scalar product, creation operators are adjoint to 
annihilation ones and vice versa, then we say that 
we have a well--defined system with generalized 
statistics.

\paragraph*{Example 1:}
We assume here that $\Pi_n \equiv P_n \equiv id_{E^{\ot n}}$.
This means that the algebra of states $\ca$ is identical with
the full tensor algebra $TE$ over the space $E$, and the second
algebra $\ca^{\ast}$ is identical with the tensor algebra 
$TE^{\ast}$. The action (\ref{act}) of annihilation operators
is given by the formula
\be
a_{x^{\ast i_k} \ot \cdots \ot x^{\ast i_1}} 
(x^{j_1} \ot \cdots \ot x^{j_n}) 
:= \d_{i_1}^{j_1} \cdots \d_{i_k}^{j_k} \
x^{j_{n-k+1}} \ot \cdots \ot x^{j_n}.
\ee
For the scalar product we have the equation
\be
\langle i_n \cdots i_1 |j_1\cdots j_n \rangle^n
:= \d^{i_1 j_1}\cdots \d^{i_n j_n}
\ee
It is easy to see that we have the relation
and
\be
a_{x^{\ast i}} a_{x_j} := \d_i^j {\bf 1}.
\ee
In this way we obtain the most simple example of 
well--defined system with generalized statistics. 
The corresponding statistics is the so--called 
infinite (Bolzman) statistics \cite{owg,gre}.

\paragraph*{Example 2:}
For this example we assume that $\Pi_n \equiv id_{E^{\ot n}}$.
This means that $\ca \equiv TE$ and $\ca^{\ast} \equiv TE^{\ast}$. 
For the scalar product and for the action of annihilation operators
we assume that there is a linear and invertible operator 
$T : E^{\ast} \ot E \ra E \ot E^{\ast}$ defined by its
matrix elements
\be
T(x^{\ast i}\ot x^j) = \Sigma_{k,\ast l} \ 
T^{\ast ij}_{k\ast l} \ x^k \ot x^{\ast l},
\label{teo}
\ee
such that we have
\be
(T^{\ast ij}_{k\ast l})^{\ast} = 
\overline{T}^{\ast ji}_{l\ast k},
\mbox{i.e.} T^{\ast} = \overline{T}^t,
\ee
and $(T^t)^{\ast ij}_{k\ast l} = T^{\ast ji}_{l\ast k}$.
Note that this operator not need to be linear, one can 
also consider the case of nonlinear one. We also 
assume that the operator $T^{\ast}$ act to the left, 
i.e. we have the relation
\be
(x^{\ast j} \ot x^{i})T^{\ast} = \Sigma_{l,\ast k} \ 
(x^{l} \ot x^{\ast k}) \ \overline{T}^{\ast ji}_{l\ast k},
\ee
and
\be
(T(x^{\ast i} \ot x^j))^{\ast} \equiv (x^{\ast j} 
\ot x^{i}) T^{\ast}.
\ee
The operator $T$ given by the formula (\ref{teo}) 
is said to be {\it a twist} or {\it a cross} operator. 
The operator $T$ describes the cross statistics
of quasiparticles and quasiholes.
The set $P$ of projectors is defined by induction
\be
P_{n+1} := (id \ot P_n) \circ R_{n+1},
\ee
where $P_1 \equiv id$ and the operator $R_n$ is given by 
the formula 
\be
R_n := id + \tilde{T}^{(1)} + \tilde{T}^{(1)} \tilde{T}^{(2)}
+ \cdots + \tilde{T}^{(1)}\dots\tilde{T}^{(n-1)} ,
\ee
where 
$\tilde{T}^{(i)} := id_E \ot \cdots 
\ot \tilde{T} \ot \cdots \ot id_E$, 
$\tilde{T}$ on the $i$--th place, and
\be
(\tilde{T})^{ij}_{kl} = T^{\ast ki}_{l\ast j}.
\ee
If the operator $\tilde{T}$ is a bounded operator acting 
on some Hilbert space such that we have the following
Yang-Baxter equation on $E\ot E \ot E$ 
\be
(\tl T, \ot id_E )\circ (id_E \ot \tl T, )\circ 
(\tl T, \ot id_E ) = (id_E \ot \tl T, ) \circ 
(\tl T, \ot id_E )\circ (id_E \ot \tl T, ),
\ee
and $||\tl T,|| \leq 1$,
then according to Bo$\dot{z}$ejko and Speicher \cite{bs2}
there is a positive definite scalar product
\be
\pg s, t, n, _T := \pg s, P_n t, n, _0
\label{csca}
\ee
for $s, t \in \ca^n \equiv E^{\ot n}$. Note that the 
existence of nontrivial kernel of operator 
$P_2 \equiv R_1 \equiv id_{E\ot E} + \tilde{T}$ 
is essential for the nondegeneracy of the scalar product
\cite{jswe}.
One can see that if this kernel is trivial, then we
obtain well--defined system with generalized statistics
\cite{m10,ral}.

\paragraph*{Example 3:}
If the kernel of $P_2$  is nontrivial, then the scalar 
product (\ref{csca}) is degenerate. Hence we must remove 
this degeneracy by factoring the mentioned above scalar 
product by the kernel. We assume that there is an ideal
$I \subset TE$ generated by a subspace 
$I_{2} \subset ker P_2 \subset E \ot E$ such that
\be
a_{s^{\ast}} I \subset I 
\ee
for every $s^{\ast} \in \ca^{\ast}$, and for the
corresponding ideal $I^{\ast} \subset E^{\ast} \ot E^{\ast}$ 
we have
\be
a_{s^{\ast}} t = 0
\ee
for every $t \in TE$ and $s^{\ast} \in I^{\ast}$.
The above ideal $I$ is said to be Wick ideal \cite{jswe}.
We have here the following formulae
\be
\ca := TE/I, \quad \ca^{\ast} := TE^{\ast}/I^{\ast}
\ee
for our algebras. The projection $\Pi$ is the quotient map
\be
\Pi : \tl s, \in TE \ra s \in TE/I \equiv \ca
\ee
For the scalar product we have here the following relation
\be
\langle s|t\rangle_{B,T} := \langle\tilde{s}|\tilde{t}\rangle_T
\ee
for $s = P_m (\tilde{s})$ and $t = P_n (\tilde{t})$.
One can define here the action of annihilation operators
in such a way that we obtain well--defined system with 
generalized statistics \cite{ral}.

\paragraph*{Example 4:}
If a linear and invertible operator $B: E \ot E \ra E \ot E$
defined by its matrix elements
\be
B(x^i \ot x^j) := B^{ij}_{kl} (x^k \ot x^l)
\ee
is given such that we have the following conditions
\be
\ba{l}
B\1n B\2n B\1n = B\2n B\1n B\2n,\\
B^{(1)}T\2nT\1n = T\2nT\1n B^{(2)},\\
(id_{E \ot E} + \tilde{T})(id_{E \ot E} - B) = 0,
\label{cd}
\ea
\ee
then one can prove that there is well defined action
of annihilation operators and scalar product. In this 
case we need two operators $T$ and $B$ satisfying the 
above consistency conditions for the model with 
generalized statistics \cite{RM,m10,ral}. 

\paragraph*{Example 5:}
If $B = \frac{1}{\mu} \tl T,$, where $\mu$ is a parameter, 
then the third condition (\ref{cd}) is equivalent to the well 
known Hecke condition for $\tl T,$ and we obtain the 
well--known relations for Hecke symmetry and quantum groups
\cite{P,PW,Ke}.
\section*{Physical applications}
Let us consider the system equipped with generalized
statistics and described by two operators $T$ and 
$B$ like in Example 4. We assume here in addition that a 
linear and Hermitian operator $S : E\ot E\ra E\ot E$ such that
\be
S\1n S\2n S\1n = S\2n S\1n S\2n, \quad \mbox{and} \quad
S^2 = id_{E\ot E}.
\ee
is given. If we have the following relation
\be
\tl T, \equiv B \equiv S,
\ee
then it is easy to see that the conditions (\ref{cd})
are satisfied and we have well--defined system with
generalized statistics.
Let us assume for simplicity that the operator $S$ is
diagonal and is given by the following equation
\be
S(x^i \ot x^j) = \ep^{ij} x^j \ot x^i,
\ee
for $i, j = 1, \ldots, N$, where $\ep^{ij} \in \com$,
and $\ep^{ij}\ep^{ji} = 1$. In the general case we have
\be
\ep_{ij} = (-1)^{\Sigma_{ij}} q^{\Om_{ij}},
\label{comf}
\ee
where $\Sigma := (\Sigma_{ij})$ and $\Om := (\Om_{ij})$ are 
integer--valued matrices such that $\Sigma_{ij} = \Sigma_{ji}$ 
and $\Om_{ij} = - \Om_{ji}$, $q \in \com \setminus \{0\}$ is
a parameter \cite{zoz}. The algebra $\ca$ is here a quadratic 
algebra generated by relations
\be
x^i \odot x^j = \ep^{ij} x^j \odot x^i,\quad\mbox{and}
\quad (x^i)^2 = 0 \quad\mbox{if}\quad\ep^{ii} = -1
\ee
We also assume that $\ep^{ii} = -1$ for every 
$i = 1, \ldots, N$. In this case the algebra $\ca$ 
is denoted by $\Lambda_{\ep}(N)$. One can see that this
is a $G$--graded $\ep$--commutative algebra \cite{sch}.
Now let us study the algebra $\Lambda_{\ep}(2)$, where 
$\ep^{ii} = -1$ for $i=1, 2$, and $\ep^{ij} = 1$ for $i\neq j$, 
in more details. In this case our algebra is generated by $x^1$ 
and $x^2$ such that we have 
\be
x^1 \odot x^2 = x^2 \odot x^1, \quad (x^1)^2 = (x^2)^2 = 0
\ee
Note that the algebra $\Lambda_{\ep}(2)$ is an example of 
the so--called $Z_2 \oplus Z_2$--graded commutative colour 
Lie superalgebra \cite{luri}. Such algebra can be transformed
into the usual grassmann algebra $\Lambda_2$ generated by 
$\th^1$ and $\th^2$ such that we have the anticommutation 
relation
\be
\th^1 \ \th^2 = - \th^2 \ \th^1,
\label{thr}
\ee
and $(\th^1)^2 = (\th^1)^2 = 0$.
In order to do such transformation we use the Clifford algebra 
$C_2$ generated by $e^1, e^2$ such that we have the relations
\be
e^i \ e^j + e^j \ e^i = 2 \d^{ij}\quad \mbox{for} \quad 
i, j = 1, 2.
\ee
For generators $x^1$, and $x^2$ of the algebra $\Lambda_{\ep}(2)$ 
the transformation is given by
\be
\th^1 := x^1 \ot e^1, \quad\mbox{and}\quad \th^2 := x^2 \ot e^2.
\ee
It is interesting that the algebra $\Lambda_{\ep}(2)$ can be 
represented by one grassmann variable $\th$, $\th^2 = 0$ 
\be
x^1 = (\th, 1),\quad  x^2 = (1, \th).
\label{sqa}
\ee
For the product $x^1 \odot x^2$ we obtain
\be
x^1 \odot x^2 = (\th, \th).
\label{haf}
\ee
In physical interpretation generators $\th^1$ and $\th^2$ of
the algebra $\Lambda_2$ represents two fermions. They anticommute
and according to the Pauli exclusion principle we can not put 
them into one energy level. Observe that the corresponding 
generators $x^1$ and $x^2$ of the algebra $\Lambda_{\ep}(2)$ 
commute, their squares disappear and they describe two different 
quasiparticles. This means that these quasiparticles behave 
partially like bosons, we can put them simultaneously into one 
energy levels. This also means that single fermion can be
transform under certain interactions into a system of two 
different quasiparticles. 


\begin{thebibliography}{}
\bibitem{owg} O. W. Greenberg, {\it Phys. Rev. Lett.} {\bf 64} 
705 (1990).
\bibitem{gre} O. W. Greenberg, {\it Phys. Rev.} {\bf D 43}, 
4111 (1991).
\bibitem{moh} R. N. Mohapatra, {\it Phys. Lett.} {\bf B 242}, 
407 (1990).
\bibitem{fi} D. I. Fivel, {\it Phys. Rev. Lett.} {\bf 65}, 
3361, (1990).
\bibitem{zag} D. Zagier, {\it Commun. Math. Phys.} {\bf 147}, 
199 (1992).
\bibitem{mepe} S. Meljanac and A. Perica, 
{\it Mod. Phys. Lett.} {\bf A9} 3293 (1994).
\bibitem{bs2} M. Bo$\dot{z}$ejko, R. Speicher,
{\it Math. Ann.} {\bf 300}, 97, (1994).
\bibitem{P} W. Pusz, {\it Rep. Math. Phys.} {\bf 27}, 
394, (1989)
\bibitem{PW} W. Pusz and S.L. Woronowicz,
{\it Rep. Math. Phys} {\bf 27}, 231, (1989)
\bibitem{ckl} M. Chaichian, P. Kulisch, J. Lukierski,
{\it Phys. Lett.} {\bf B262}, 43, (1991).
\bibitem{V} S.P. Vokos, {\it J. Math. Phys.} {\bf 32}, 2979, (1991).
\bibitem{FZ} D.B. Fairle and C.K. Zachos,
{\it Phys. Lett.} {\bf B256}, 43, (1991)
\bibitem{Wu}  Y.S. Wu, {\it J.Math.Phys.} {\bf 52}, 2103, 1984
\bibitem{I} T.D. Imbo and J. March--Russel,
{\it Phys. Lett.} {\bf B252}, 84, 1990
\bibitem{SM} S. Majid, {\it Int. J. Mod. Phys.}{\bf A5}, 1 (1990).
\bibitem{Maj} S. Majid, {\it J. Math. Phys.}{\bf 34}, 1176, (1993).
\bibitem{SMa} S. Majid, {\it J. Math. Phys.}{\bf 34}, 4843, (1993).
\bibitem{bm} S. Majid, {\it J. Math. Phys.}{\bf 34}, 2045 (1993).
\bibitem{sma} S. Majid, Algebras and Hopf Algebras in Braided
Categories, in Advanced in Hopf Algebras, Plenum 1993.
\bibitem{qm} S. Majid, {\it J. Geom. Phys.} {\bf 13}, 169 (1994).
\bibitem{ssm} S. Majid, {\it AMS Cont. Math.} {\bf 134}, 219 (1992).
\bibitem{WM3} W. Marcinek, {\it J. Math. Phys.} {\bf 33}, 
1631 (1992).
\bibitem{WM4} W. Marcinek, {\it Rep. Math. Phys.} {\bf 34}, 
325 (1994).
\bibitem{WM6} W. Marcinek, {\it Rep. Math. Phys.} {\it 33}, 
117, (1993).
\bibitem{WM8} W. Marcinek, {\it J. Math. Phys.} {\bf 35}, 
2633, (1994).
\bibitem{wm7} W. Marcinek, {\it Int. J. Mod. Phys.}
{\bf A10}, 1465-1481 (1995).
\bibitem{mad} W. Marcinek, On the deformation of commutation
relations, in Proceedings of the XIII Workshop in Geometric 
Methods in Physics, July 1-7, 1994 Bia{\l}owieza, Poland, 
ed. J. Antoine, Plenum Press 1995.
\bibitem{mco} W. Marcinek, On algebraic model of composite 
fermions and bosons, in  Proceedings of the IXth Max Born 
Symposium, Karpacz, September 25 - September 28, 1996, Poland. 
\bibitem{qweyl} W. Marcinek, On quantum Weyl algebras and 
generalized quons, in Proceedings of the symposium:  Quantum 
Groups and Quantum Spaces, Warsaw, November 20-29, 1995, 
Poland, ed. by R. Budzynski, W. Pusz and S. Zakrzewski, 
Banach Center Publications, Warsaw 1997.
\bibitem{wmq} W. Marcinek, {\it Rep. Math. Phys.} {\bf 41}, 
155 (1998).
\bibitem{jswe} P.E.T. Jorgensen, L.M. Schmith, and
R.F. Werner, Positive representation of general
commutation relations allowing Wick ordering, 
{\it J. Funct. Anal.} {\bf 134}, 33 (1995).
\bibitem{RM} W. Marcinek and R. Ra{\l }owski, Particle operators
from braided geometry, in "Quantum Groups, Formalism
and Applications"  XXX Karpacz Winter School in Theoretical 
Physics, 1994, Eds. J. Lukierski et al., 149-154 (1995).
\bibitem{m10} W. Marcinek and Robert Ra{\l }owski, On Wick
Algebras with Braid Relations, {\it Preprint IFT UWr} 
{\bf 876/9}, (1994)
and {\it J. Math. Phys.} {\bf 36}, 2803, (1995).
\bibitem{ral} R. Ralowski, {\it J. Phys.}{\bf A30}, 2633 (1997).
\bibitem{sci} R. Scipioni, {\it Phys. Lett.} {\bf B327}, 56 (1994).
\bibitem{twy} Yu Ting and Wu Zhao-Yan, {\it Science in China}
{\bf A37}, 1472 (1994).
\bibitem{mep} S. Meljanac and A. Perica {\it Mod. Phys. Lett.}
{\bf A9}, 3293 (1994).
\bibitem{mphi} M. Pillin, {\it Commun. Math. Phys.} {\bf 180},
23 (1996).
\bibitem{mira} A. K. Mishra and G. Rajasekaran,
{\it J. Math. Phys.} {\bf 38}, 466 (1997).
\bibitem{fios} G. Fiore and P. Schup, Statistics and Quantum
Group Symmetries, in Proceedings of the symposium:  Quantum 
Groups and Quantum Spaces, Warsaw, November 20-29, 1995, 
Poland, ed. by R. Budzynski, W. Pusz and S. Zakrzewski, 
Banach Center Publications, Warsaw 1997.
\bibitem{melme} S. Meljanac and M. Molekovic, 
{\it Int. J. Mod. Phys. Lett.} {\bf A11}, 139 (1996).
\bibitem{castat} W. Marcinek, Categories and quantum statistics,
in Proceedings of the symposium: Quantum Groups and their 
Applications in Physics, Poznan October 17-20, 1995, Poland, 
{\it Rep. Math. Phys.} {\bf 38}, 149-179 (1996)
\bibitem{top} W. Marcinek, Topology and quantization, in 
Proceedings of the IVth International School on Theoretical 
Physics, Symmetry and Structural Properties, Zajaczkowo 
k. Poznania, August 29 - September 4 1996, Poland.
\bibitem{sin} W. Marcinek, {\it J. Math. Phys.}
{\bf 39}, 818--830 (1998).
\bibitem{zee} A. Zee, Quantum Hall fluids in Field Theory,
Topology and Condensed Matter Physics, ed. by H. D. Geyer,
Lecture Notes in Physics, Springer 1995.
\bibitem{jai} J. K. Jain, {\it Phys. Rev. Lett.} {\bf 63}, 
199 (1989),
{\it Phys. Rev.} {\bf B 40}, 8079 (1989); {\bf 41}, 7653 (1990).
\bibitem{dst} R. R. Du, H. L. St\"ormer, D. C. Tsui, A. S. Yeh,
L. N. Pfeiffer and K. W. West, {\it Phys. Rev. Lett.} 
{\bf 73}, 3274 (1994).
\bibitem{hal1} F. D. M. Haldane, {\it J. Phys.}
{\bf C14}, 2585 (1981).
\bibitem{bys} K. Byczuk and J. Spalek, 
{\it Phys. Rev.} {\bf B51}, 7934 (1995).
\bibitem{Ke} A. Kempf, {\it Let. Math. Phys.} {\bf 26}, 
11, (1992).
\bibitem{zoz} Z. Oziewicz, Lie algebras for arbitrary grading group, in
Differential Geometry and Its Applications ed. by J. Janyska and D. Krupka,
World Scientific, Singapore 1990
\bibitem{sch} M. Scheunert, {\em J. Math. Phys.} {\bf 20}, 712, (1979)
\bibitem{luri} J. Lukierski, V. Rittenberg, 
{\it Phys. Rev.} {\bf D18}, 385, (1978).
\end{thebibliography}
\end{document}